\documentclass[11pt]{amsart}
\usepackage{amssymb,amsmath,amsgen,amsxtra,dsfont,times,a4wide} 
\begin{document}
%
%
%
\theoremstyle{definition}
\newtheorem{Definition}{Definition}[section]
\newtheorem{Example}[Definition]{Example}
\newtheorem{Examples}[Definition]{Examples}
\newtheorem{Remark}[Definition]{Remark}
\newtheorem{Remarks}[Definition]{Remarks}
\newtheorem{Caution}[Definition]{Caution}
\newtheorem{Conjecture}[Definition]{Conjecture}
\newtheorem{Question}[Definition]{Question}
\newtheorem{Questions}[Definition]{Questions}
\theoremstyle{plain}
\newtheorem{Theorem}[Definition]{Theorem}
\newtheorem{Proposition}[Definition]{Proposition}
\newtheorem{Lemma}[Definition]{Lemma}
\newtheorem{Corollary}[Definition]{Corollary}
\newtheorem{Fact}[Definition]{Fact}
\newtheorem{Facts}[Definition]{Facts}
\newtheoremstyle{voiditstyle}{3pt}{3pt}{\itshape}{\parindent}%
{\bfseries}{.}{ }{\thmnote{#3}}%
\theoremstyle{voiditstyle}
\newtheorem*{VoidItalic}{}
\newtheoremstyle{voidromstyle}{3pt}{3pt}{\rm}{\parindent}%
{\bfseries}{.}{ }{\thmnote{#3}}%
\theoremstyle{voidromstyle}
\newtheorem*{VoidRoman}{}
%
\newcommand{\prf}{\par\noindent{\sc Proof.}\quad}
\newcommand{\blowup}{\rule[-3mm]{0mm}{0mm}}
\newcommand{\Aff}{{\mathds{A}}}
\newcommand{\BB}{{\mathds{B}}}
\newcommand{\CC}{{\mathds{C}}}
\newcommand{\FF}{{\mathds{F}}}
\newcommand{\GG}{{\mathds{G}}}
\newcommand{\HH}{{\mathds{H}}}
\newcommand{\NN}{{\mathds{N}}}
\newcommand{\ZZ}{{\mathds{Z}}}
\newcommand{\PP}{{\mathds{P}}}
\newcommand{\QQ}{{\mathds{Q}}}
\newcommand{\RR}{{\mathds{R}}}
\newcommand{\Sphere}{{\mathds{S}}}
\newcommand{\lin}{\sim}
\newcommand{\num}{\equiv}
\newcommand{\dual}{\ast}
\newcommand{\iso}{\cong}
\newcommand{\caniso}{\cong}
\newcommand{\mm}{{\mathfrak m}}
\newcommand{\pp}{{\mathfrak p}}
\newcommand{\qq}{{\mathfrak q}}
\newcommand{\rr}{{\mathfrak r}}
\newcommand{\pP}{{\mathfrak P}}
\newcommand{\qQ}{{\mathfrak Q}}
\newcommand{\rR}{{\mathfrak R}}
%
%
\newcommand{\perdef}{{\stackrel{{\rm def}}{=}}}
\newcommand{\Sym}{{\mathfrak S}}
\newcommand{\EnSym}{{\cal S}}
\newcommand{\Cycl}[1]{{\ZZ/{#1}\ZZ}}
\newcommand{\Elem}[1]{{{\rm E}_{#1}}}
\newcommand{\myExt}{{\mathcal E}}
\newcommand{\myKernel}{{\mathcal K}}
\newcommand{\myBigExt}{\widetilde{{\mathcal E}}}
\newcommand{\myBigKernel}{\widetilde{{\mathcal K}}}
\newcommand{\Aut}{{\rm Aut}}
\newcommand{\Hom}{{\rm Hom}}
\newcommand{\ab}{{\rm ab}}
\newcommand{\Xaff}{{X^{\rm aff}}}
\newcommand{\fgal}{{f_{\rm gal}}}
\newcommand{\Xgal}{{X_{\rm gal}}}
\newcommand{\Xgalaff}{{X_{\rm gal}^{\rm aff}}}
\newcommand{\Xaffgal}{\Xgalaff}
\title[Natural Central Extensions]{Natural Central Extensions of Groups}
\author{Christian Liedtke}
\address{Mathematisches Institut, Heinrich-Heine-Universit\"at,
  40225 D\"usseldorf, Germany}
\email{liedtke@math.uni-duesseldorf.de}
\date{September 13, 2007}
\begin{abstract}
  Given a group $G$ and an integer $n\geq2$ we construct a new group
  $\myBigKernel(G,n)$.
  Although this construction naturally occurs in the context of finding new
  invariants for complex algebraic surfaces,
  it is related to the theory of central extensions and the Schur multiplier.
  A surprising application is that Abelian groups of odd order possess
  naturally defined covers that can be computed from a given cover by a kind 
  of warped Baer sum.
\end{abstract}
\maketitle
\tableofcontents
\section*{Introduction}

Given a group $G$ and an integer $n\geq2$ we 
introduce a new group $\myBigKernel(G,n)$.
\medskip

It originates from the author's work \cite{li} on 
Moishezon's programme \cite{mo} for complex algebraic surfaces.
More precisely, to obtain finer invariants for surfaces,
one attaches to a surface and an embedding into projective space
the fundamental group of $\pi_1(\PP^2-D)$, where
$D$ is a curve (the branch locus of a generic projection)
in the projective plane $\PP^2$.
Although knowing this fundamental group (and its monodromy morphism, see
Section \ref{surfacesection} for a precise statement) one can reconstruct
the given surface, this fundamental group is too complicated and too
large to be useful.

Instead one looks for subgroups and subquotients of this group
$\pi_1(\PP^2-D)$ to obtain the desired invariants.
The most prominent one is a naturally defined subquotient that
has itself a geometric interpretation, namely the fundamental
group of the Galois closure of a generic projection from
the given surface.
Its computation already in special cases by Moishezon and Teicher
\cite{mote} led to counter-examples to Bogomolov's watershed
conjecture.
\medskip

Our construction $\myBigKernel(G,n)$ is closely related to the fundamental
group of the Galois closure just mentioned 
(see Section \ref{surfacesection} for an exposition).
Here, we will be concerned with its group theoretical properties:

In general, it is difficult to compute $\myBigKernel(G,n)$ for given
$G$ and $n$.
For example, if $G$ is Abelian we usually obtain a nilpotent group
of class $2$.
On the other hand, $\myBigKernel(G,n)$ inherits many properties from $G$:
for example, if $G$ is finite, nilpotent, or solvable, then the same
will be true for $\myBigKernel(G,n)$ and for all $n\geq2$.
This construction, defined in Section \ref{mainsection}, is
in the spirit of Hopf's computation of $H_2(G,\,\ZZ)$ in terms of 
a presentation of $G$.
And so it is not surprising there is a connection with the
theory of central extensions, covers and the Schur multiplier.

Let $G$ be a finite group and $M:=M(G)$ be its Schur multiplier.
Then a cover of a finite group $G$ is defined to be a group 
$G^\ast$, which is a central extension
$$\begin{array}{ccccccccc}
  0&\to&M&\to&G^\ast&\to&G&\to&1
\end{array}$$
of $M$ by $G$ such that $M$ lies in the commutator 
subgroup of $G^\ast$.
The notion of covers comes from the study of central 
extensions and projective representations, see e.g. \cite{ka}.

Every perfect group has unique cover, 
which is called its universal central extension.
However, non-perfect groups usually do not have unique or natural 
covers.
\medskip

If we apply our construction to an Abelian group $G$, we obtain
for every $n\geq2$ a naturally defined central extension
$$\begin{array}{ccccccccc}
  0&\to&M(G)&\to&\myBigKernel(G,n)&\to&G^{n-1}&\to&1\,.
\end{array}$$
If $n\geq3$ or if $n=2$ and $G$ has odd order, then $M(G)$ lies in the
centre of $\myBigKernel(G,n)$.
\medskip

A surprising application (the case $n=2$) is the existence of
natural covers for Abelian groups of odd order.
This natural cover can be computed from a given cover by a sort 
of warped Baer sum of this given cover with itself.
\medskip

This article is organised as follows.

In Section \ref{elementarysection} we construct the
auxiliary group $\myKernel(G,n)$ for $n\geq2$ and a given group $G$.
It is a subgroup of $G^n$ and not so difficult to compute, 
especially when $G$ is perfect or Abelian.

In Section \ref{mainsection} we use this auxiliary construction
to define $\myBigKernel(G,n)$, the main object of this article.
We show that it is a central extension of $H_2(G,\ZZ)$ 
by $\myKernel(G,n)$.
Also, we prove that it inherits properties such as finiteness,
nilpotency or solvability from $G$.

In Section \ref{centralsection} we prove that 
$\myBigKernel(G,n)$ can be computed from an arbitrary cover of $G$.
In particular, it contains the universal central extension of $G$
in case $G$ is perfect.

In Section \ref{abeliansection}
we determine the structure, i.e. the centre, Frattini and 
commutator subgroup, of $\myBigKernel(G,n)$ in case $G$ is Abelian.
Here we also prove that Abelian groups of odd order possess
natural covers.

In Section \ref{surfacesection} we describe the relation to
fundamental groups groups of algebraic surfaces
and to Moishezon's programme to find finer invariants for
surfaces.

\begin{VoidRoman}[Acknowledgements]
This article extends results from my Ph.D. thesis \cite{li}.
I thank my supervisor Gerd~Faltings for discussions and
help, as well as
the Max-Planck-Institut in Bonn for hospitality and financial support.
Also, I thank Irene~Bouw for proof-reading and many suggestions.
\end{VoidRoman}

\section{An auxiliary construction}
\label{elementarysection}

We let $G$ be a group and $n\geq2$ a natural number.
We denote by $G^\ab:=G/\left[G,G\right]$ the Abelianisation of $G$.
Then we define a map
$$\begin{array}{ccccc}
        \psi&:&G^n&\to&G^\ab \\
        &&(g_1,...,g_n)&\mapsto&\overline{g_1\cdot...\cdot g_n}
  \end{array}$$
which is a homomorphism as $G^\ab$ is Abelian.
\begin{Definition}
  For a group $G$ and a natural number $n\geq2$ we define
  $\myKernel(G,\,n)$ to be the kernel of the homomorphism
  $\psi:G^n\to G^\ab$.
\end{Definition}
Clearly $\myKernel(-,n)$ is functorial in its first argument.
We start with the
\begin{Proposition}
  \label{mykernelproperties}
  Let $G_1$, $G_2$, $G$ be arbitrary groups and $n\geq2$ a natural number.
  \begin{enumerate}
  \item If $G_1\to G_2$ is an injective or a surjective homomorphism then 
    the same is true for the induced maps $\myKernel(G_1,n)\to\myKernel(G_2,n)$.
  \item There exists a natural isomorphism
  $\myKernel(G_1\times G_2,n)\caniso\myKernel(G_1,n)\times\myKernel(G_2,n)$.
  \item For $n\geq3$ the natural homomorphism from $\myKernel(G,n)^\ab$
    onto $\myKernel(G^\ab,n)$ is an isomorphism.
  \end{enumerate}
\end{Proposition}

\prf
The first two assertions follow immediately from the definition.

The surjection $G\to G^\ab$ and the universal property of the Abelianisation
imply that there is a natural surjective homomorphism
$\myKernel(G,n)^\ab\to\myKernel(G^\ab,n)$.
An element of the kernel $\myKernel(G,n)\to\myKernel(G^\ab,n)$ is also
an element of the kernel $G^n\to(G^\ab)^n$, which is $[G,G]^n$.
Since we assumed that $n\geq3$, we may write
\begin{equation}
\label{commutatoreq}
\begin{array}{cccc}
  ([h_1,h_2],1,...,1) &=& [ (h_1,h_1{}^{-1},1,...,1),\,(h_2,1,h_2{}^{-1},...,1) ] 
  &\in G^n.
\end{array}
\end{equation}
Thus $[G,G]^n$ is not only a subgroup of $\myKernel(G,n)$ but also lies
inside the commutator subgroup of $\myKernel(G,n)$.
Hence the kernel $\myKernel(G,n)\to\myKernel(G^\ab,n)$ is the commutator
subgroup of $\myKernel(G,n)$ and we are done.
\qed

\begin{Remark}
  Already here we see that the case $n=2$ has to be treated separatedly.
  If we need $n\geq3$ for a statement it is usually easy to obtain a
  counter-example for the corresponding statement for $n=2$ 
  using elementary Abelian $2$-groups, dihedral groups or 
  the quaternion group.
\end{Remark}

In the following two cases it is particularly easy to determine 
$\myKernel(G,n)$.
\begin{Proposition}
  \label{mykernelAbelian}
  \begin{enumerate}
  \item 
    If $G$ is perfect then $\myKernel(G,n)\caniso G^n$.
  \item
    If $G$ is Abelian then 
    $\myKernel(G,n)\iso G^{n-1}$.
    This isomorphism is not canonical.
  \end{enumerate}
\end{Proposition}
\prf
The Abelianisation of a perfect group is trivial and so the
first assertion follows from the definition of $\myKernel(G,n)$.

Now let $G$ be an Abelian group.
Then the map
$$\begin{array}{ccc}
  G^{n-1}&\to& G^n\\
  (g_1,...,g_{n-1})&\mapsto&(g_1,...,g_{n-1},(g_1\cdot...\cdot g_{n-1})^{-1})
\end{array}$$
defines a homomorphism.
It is injective with image $\myKernel(G,n)$.
\qed

\begin{Proposition}
  \label{mykernelheritage}
  Let $n\geq2$ and let $P$ be one of the following properties:
  \begin{center}
    Abelian,\quad finite,\quad nilpotent,\quad perfect,\quad solvable.
  \end{center}
  Then $G$ has the property $P$ if and only if $\myKernel(G,n)$
  has the same property.
\end{Proposition}
\prf
By definition, $\myKernel(G,n)$ is a subgroup of $G^n$.
Therefore if $G$ is Abelian (resp. finite, nilpotent, solvable) the
same is true for $\myKernel(G,n)$.
If $G$ is perfect then $\myKernel(G,n)=G^n$, which is also perfect.

The projection onto the first factor $G^n\to G$ induces a surjective
homomorphism from $\myKernel(G,n)$ onto $G$.
Hence $G$ is a quotient of $\myKernel(G,n)$.
Therefore if $\myKernel(G,n)$ is Abelian (resp. finite, nilpotent, perfect,
solvable) the same is true for $G$.
\qed

\section{The main construction}
\label{mainsection}

As in the previous section, we let $G$ be a group and $n\geq2$ be a 
natural number.
We choose a presentation $G\iso F/N$ where $F$ is a free group.
Then $\myKernel(N,n)$ is a subgroup of $\myKernel(F,n)$ which is a
subgroup of $F^n$.

We denote by $\ll\myKernel(N,n)\gg$ the subgroup normally generated
by $\myKernel(N,n)$ inside $F^n$.
For $n\geq3$, it is not difficult to see (using formula (\ref{commutatoreq}))
that the normal closure $\ll\myKernel(N,n)\gg$ of $\myKernel(N,n)$
inside $F^n$ is equal to the normal closure of $\myKernel(N,n)$ inside
$\myKernel(F,n)$.

\begin{Definition}
  We let $G$ be a group and $n\geq2$ be a natural number.
  We define
  $$\begin{array}{ccc}
    \myBigKernel(G,n) &:=& \myKernel(F,n)/\ll\myKernel(N,n)\gg\,.
  \end{array}$$
\end{Definition}

\begin{Theorem}
  \label{mybigkernelthm}
  The group $\myBigKernel(G,n)$ does not depend upon the choice
  of a presentation.
  There exists a central short exact sequence
  \begin{equation}
  \label{mybigkerneleq}
  \begin{array}{cccccccccc}
    0&\to&H_2(G,\,\ZZ)&\to&\myBigKernel(G,\,n)&\to&\myKernel(G,n)&\to&1&.
  \end{array}
  \end{equation}
  For $n\geq3$ the group $H_2(G,\ZZ)$ lies inside the commutator
  subgroup of $\myBigKernel(G,n)$.
\end{Theorem}
\prf
We choose a presentation $G\iso F/N$ and abbreviate
the normal closure $\ll\myKernel(N,n)\gg$ of $\myKernel(N,n)$ 
in $F^n$ by $R$.

First, we will prove the short exact sequence of the statement
of the theorem:
Let $\pi$ be the projection of $F^n$ onto its last
$n-1$ factors. 
By abuse of notation we denote its restriction to $\myKernel(F,n)$ 
again by $\pi$.
We obtain a short exact sequence
$$\begin{array}{ccccccccc}
  1&\to&[F,F]&\to&\myKernel(F,n)&\stackrel{\pi}{\to}&F^{n-1}&\to&1.
\end{array}$$
An easy computation with commutators shows
that $R\cap\ker\pi=[F,N]$.
Via $\pi$ we obtain the following diagram with exact rows and columns:
\begin{equation}
\label{mainprfeq}
\begin{array}{ccccccccc}
  1&\to&[F,N]&\to&R&\to&N^{n-1}&\to&1 \\
  &&\downarrow&&\downarrow&&\downarrow\\
  1&\to&N\cap[F,F]&\to&N^n\cap\myKernel(F,n)&\to&N^{n-1}&\to&1 \\
  &&\downarrow&&\downarrow&&\downarrow\\
  1&\to&[F,F]&\to&\myKernel(F,n)&\stackrel{\pi}{\to}&F^{n-1}&\to&1 
\end{array}
\end{equation}
Taking quotients of successive rows we exhibit
$\myKernel(F,n)/R$ as an
extension of $(N\cap[F,F])/[F,N]$ by $\myKernel(F,n)/(N^n\cap\myKernel(F,n))$.
The latter group is isomorphic to $\myKernel(G,n)$.
By Hopf's theorem (cf. \cite[Theorem II.5.3]{br}), the group
$(N\cap[F,F])/[F,N]$ is isomorphic to $H_2(G,\ZZ)$.
Hence we obtain an extension
$$\begin{array}{ccccccccc}
  1&\to&H_2(G,\,\ZZ)&\to&\myKernel(F,n)/R&\to&\myKernel(G,n)&\to&1.
\end{array}$$

Next, we will show that this extension is central:
Every element of $H_2(G,\ZZ)$ can be lifted to an element
of the form $\vec{x}:=(x,1,...,1)$ of $\myKernel(F,n)$ 
with $x\in N\cap[F,F]$.
For $\vec{y}:=(y_1,...,y_n)\in\myKernel(F,n)$ we compute
$$
\vec{y}\vec{x}\vec{y}^{-1}\,=\,
(\underbrace{[y_1,x]}_{\in[F,N]},1,...,1)\cdot(x,1,...,1) \,\equiv\,
\vec{x}\mod[F,N].
$$
Hence $H_2(G,\ZZ)$ lies inside the centre of $\myKernel(F,n)/R$.

We now prove that $\myBigKernel(G,n)$ is well-defined:
Let $\alpha:F/N\iso F'/N'$ be another presentation of $G$.
We lift this isomorphism to a map $\varphi:F\to F'$.
Then $\varphi$ maps $N$ to $N'$ and hence $\myKernel(N,n)$ to
$\myKernel(N',n)$.
Let $R'$ be the normal closure of $\myKernel(N',n)$ inside ${F'}^n$.
Then $\varphi$ induces a homomorphism
$$\begin{array}{ccccc}
\overline{\varphi} &:&\myKernel(F,n)/R &\to&\myKernel(F',n)/R'\,.
\end{array}$$
We let $\varphi'$ be another map lifting $\alpha$ to a homomorphism
from $F$ to $F'$.

Suppose now that $n\geq3$.
Then elements of the form $(1,...,f,1,...,f^{-1},1,...)$ 
generate $\myKernel(F,n)$
and so in this case it is enough to compare the maps induced by
$\varphi$ and $\varphi'$ on such elements.
For $f\in F$ there exists $n_f'\in N'$ such that
$\varphi(f)=\varphi'(f) n_f'$.
Hence
$$\begin{array}{lcl}
  \varphi((f,\,f^{-1},\,...)) &=& 
  (\varphi'(f)n_f',\, n_f'^{-1}\varphi'(f)^{-1},\,...) \\
  &=&\varphi'((f,\,f^{-1},\,...)) \cdot
  \underbrace{(n_f',\,\varphi'(f)n_f'^{-1}\varphi'(f)^{-1},\,...)}_{\in R'}.
\end{array}$$
Hence the induced maps coincide.
For $n=2$, the group $\myKernel(F,2)$ is generated by elements of the
form $(f,f^{-1})$ and $(\left[f_1,f_2\right],1)$.
It is easy to see that also in this case the induced
maps coincide.

In particular, if we choose $F=F'$ and $N=N'$ with $\alpha$ and
$\varphi$ the identity then every other lift $\varphi'$ of the
identity induces the identity on $\myKernel(F,n)/R$.

Coming back to the general case, we let
$F/N$ and $F'/N'$ be two presentations of
$G$ and let $\alpha$ be an isomorphism between them.
Then $\alpha$ and $\alpha^{-1}$ induce maps between
$\myKernel(F,n)/R$ and $\myKernel(F',n)/R'$ such that
the composites of these induced maps have to be the 
identity by the previous paragraph.
Hence $\alpha$ induces an isomorphism from
$\myKernel(F,n)/R$ to $\myKernel(F',n)/R'$.
Thus, $\myBigKernel(G,n)$ is well-defined.

Taking the quotient of the top row by the bottom row of 
(\ref{mainprfeq}) we obtain a short exact sequence
\begin{equation}
\label{perfecteq}
\begin{array}{ccccccccc}
  1&\to&[F,F]/[F,N]&\to&\myKernel(F,n)/R&\to&G^{n-1}&\to&1.
\end{array}
\end{equation}
The inclusion of $H_2(G,\ZZ)$ into $\myKernel(F,n)/R$ factors over 
$[F,F]/[F,N]$.
Suppose now that $n\geq3$.
Then $[F,F]$ lies inside the commutator subgroup of 
$\myKernel(F,n)$, cf. formula (\ref{commutatoreq}).
Hence $H_2(G,\ZZ)$ lies inside the commutator subgroup of $\myKernel(F,n)/R$.
\qed

\begin{Corollary}
  A homomorphism $\alpha:G\to H$ induces a map
  $\myBigKernel(G,n)\to\myBigKernel(H,n)$.
  
  The short exact sequence (\ref{mybigkerneleq}) 
  induces maps
  $H_2(G,\ZZ)\to H_2(H,\ZZ)$ and 
  $\myKernel(G,n)\to\myKernel(H,n)$.
  These maps coincide with the map induced by $\alpha$ on homology 
  and the map induced by $\alpha$ from $\myKernel(G,n)$ to $\myKernel(H,n)$,
  respectively.
\end{Corollary}
\prf
We choose presentations $G\iso F/N$ and $H\iso F'/N'$.
In the proof of Theorem \ref{mybigkernelthm} we did not need 
that the map $\alpha$ considered there was an
isomorphism to prove that it induces a unique map from 
$\myBigKernel(G,n)$ to $\myBigKernel(H,n)$.
This shows functoriality.

It is easy to see that the induced map coming from 
$\myBigKernel(-,n)$ is compatible with
the map induced by $\alpha$ from $\myKernel(G,n)$ to $\myKernel(H,n)$.

We have to prove that the homomorphism induced on homology
is compatible with the one coming from $\myBigKernel(-,n)$.
However, this follows from \cite[Exercise II.6.3.b]{br}.
\qed

\begin{Corollary}
  \label{bigabelianise}
  For $n\geq3$ there exist isomorphisms
  $$
  \myBigKernel(G,\,n)^\ab\,\caniso\,\myKernel(G,\,n)^\ab\,\caniso\,
  \myKernel(G^\ab,\,n)\,\iso\,(G^\ab)^{n-1}.
  $$
\end{Corollary}
\prf
The first isomorphism follows from the fact that
$H_2(G,\ZZ)$ lies inside the commutator subgroup of
$\myBigKernel(G,n)$.
The remaining isomorphisms follow from Proposition \ref{mykernelproperties}
and Proposition \ref{mykernelAbelian}.
\qed

\begin{Corollary}
  \label{mybigkernelAbelian}
  If $G$ is cyclic then 
  $\myBigKernel(G,n)\iso G^{n-1}$.
  This isomorphism is not canonical.
\end{Corollary}
\prf
If $G$ is cyclic then $H_2(G,\ZZ)$ vanishes.
Hence $\myBigKernel(G,n)$ is isomorphic to $\myKernel(G,n)$, which
is isomorphic to $G^{n-1}$ by Proposition \ref{mykernelAbelian}.
\qed

\begin{Proposition}
  Let $P$ be one of the following properties:
  \begin{center}
    finite,\quad nilpotent,\quad perfect,\quad solvable.
  \end{center}
  Then $G$ has the property $P$ if and only if 
  $\myBigKernel(G,n)$
  has the same property.
\end{Proposition}
\prf
If $G$ is finite then so are
$H_2(G,\ZZ)$ and $\myKernel(G,n)$.
Hence $\myBigKernel(G,n)$ is finite
because it is an extension of $H_2(G,\ZZ)$ by
$\myKernel(G,n)$.

Since $H_2(G,\ZZ)$ is Abelian it is nilpotent.
Hence if $G$ is solvable (resp. nilpotent) then so is
$\myBigKernel(G,n)$ because it is a (central) extension of 
two solvable (resp. nilpotent) groups.

If $G$ is perfect and $G\iso F/N$ then also 
$\tilde{G}:=\left[F,F\right]/\left[F,N\right]$ is perfect.
By the short exact sequence (\ref{perfecteq}) the group
$\myBigKernel(G,n)$ is an extension of $\tilde{G}$ by $G^{n-1}$.
Thus $\myBigKernel(G,n)$ is perfect being an extension
of two perfect groups.

The group $G$ is a quotient of $\myBigKernel(G,n)$.
So, if $\myBigKernel(G,n)$ is finite (resp. nilpotent, perfect,
solvable) the same is true for $G$.
\qed\medskip

We end this section by a remark on group actions on $\myKernel(-,n)$
and $\myBigKernel(-,n)$.

Given a group $F$, the symmetric group $\Sym_n$ on $n$ letters acts on 
$F^n$ by permuting its $n$ factors.
Clearly, this action preserves $\myKernel(F,n)$.
It is not difficult to see that if $G\iso F/N$ is a presentation of
$G$, then the $\Sym_n$-action on $\myKernel(F,n)$ induces a
$\Sym_n$-action on $\myBigKernel(G,n)$ that does not depend on the
choice of a presentation of $G$.

We let $\Sym_{n-1}$ be the subgroup of $\Sym_n$ of those permutations that
fix, say, the first letter.
Inside $\myKernel(G,n)$ (resp. $\myBigKernel(G,n)$) we form 
the normal closure $N$ (resp. $\tilde{N}$)
of the subgroup generated by the elements
$g\cdot\sigma(g^{-1})$, for all $\sigma\in\Sym_{n-1}$ and all
$g\in\myKernel(G,n)$ (resp. $g\in\myBigKernel(G,n)$).
Then the quotients $\myKernel(G,n)/N$ and $\myBigKernel(G,n)/\tilde{N}$
are isomorphic to $G$.

Thus, the $\Sym_n$-actions on $\myKernel(G,n)$ and $\myBigKernel(G,n)$
allow us to recover $G$ as a quotient of these groups.
Although we do not need this result here, it is crucial in the context of
the geometric origin of these groups.
We refer to \cite[Section 5]{li} for proofs.

\section{Central extensions and covers}
\label{centralsection}

We recall that a group $G^\ast$ is called 
a {\em cover} (or a {\em representation group}) of the 
finite group $G$ if there exists a central short exact sequence
$$\begin{array}{ccccccccc}
  0&\to&M&\to&G^\ast&\to&G&\to&1
\end{array}$$
with $M\leq \left[G^\ast,G^\ast\right]$ and such that $M$ is
isomorphic to the Schur multiplier of $G$.
For a perfect group there exists a unique cover up to isomorphism,
which is called its {\em universal central extension}.

If $G$ is finite then Pontryagin duality provides us
with a non-canonical isomorphism of its Schur multiplier 
$M(G):=H^2(G,\CC^\ast)$ with $H_2(G,\ZZ)$.

\begin{Proposition}
  \label{bigperfect}
  If $G$ is a finite and perfect group there exists a 
  short exact sequence
  $$\begin{array}{ccccccccc}
    1&\to&\widetilde{G}&\to&\myBigKernel(G,\,n)&\to&G^{n-1}&\to&1 \,.
  \end{array}$$
  Here, $\widetilde{G}$ denotes the universal central extension of $G$.
\end{Proposition}\maketitle
\prf
If $G$ is perfect with presentation $F/N$ then its universal
central extension is isomorphic to 
$\left[F,F\right]/\left[F,N\right]$, cf. \cite[Theorem 2.10.3]{ka}.
The statement follows from the short exact sequence
(\ref{perfecteq}).
\qed\medskip

This result suggests that there is a connection of $\myBigKernel(-,n)$
with the theory of central extensions.
This is in fact true by the following theorem which tells us that
we can compute $\myBigKernel(G,n)$ using an arbitrary cover of $G$.

\begin{Theorem}
  \label{coverthm}
  Let $G^\ast$ be a cover of the finite group
  $G$ and $M$ be the
  kernel of the map from $G^\ast$ onto $G$.
  For $n\geq2$ there exists an isomorphism
  $$\begin{array}{ccc}
    \myKernel(G^\ast,\,n)/\myKernel(M,\,n) &\iso&
    \myBigKernel(G,\,n).
  \end{array}$$
  In particular, the group on the left depends on $G$ and $n$
  only.
\end{Theorem}

\prf
By Schur's theorem \cite[Theorem 2.4.6]{ka},
there exists a free group $F$ and two
normal subgroups $N$ and $S$ such that
$G\iso F/N$, $G^\ast\iso F/S$ and
$N/\left[F,N\right]=S/\left[F,N\right]\times
(N\cap \left[F,F\right])/\left[F,N\right]$.

First, we show that
$$\begin{array}{ccc}
H &:=& \langle \myKernel(F,n)\cap S^n,\,\myKernel(N,n) \rangle
\end{array}$$
is a normal subgroup of $F^n$ contained in
$\myKernel(F,n)$.
Since both $S^n$ and $\myKernel(F,n)$ are normal in $F^n$ we
see that $\myKernel(F,n)\cap S^n$ is a normal subgroup of
$F^n$ contained in $\myKernel(F,n)$.
This already shows that $H$ is a subgroup of $\myKernel(F,n)$.
To show normality in $F^n$ it is enough to show that conjugates
of $\myKernel(N,n)$ by elements of $F^n$ lie inside $H$.
So let $(x_1,...,x_n)\in\myKernel(N,n)$ and $f\in F$.
Then
$$
(f,1,...)\,(x_1,...,x_n)\,(f,1,...)^{-1}\,\,=\,
\underbrace{([f,x_1],1,...)}_{\in [F,N]}\,\cdot\,\underbrace{(x_1,...,x_n)}_{\in\myKernel(N,n)}\,.
$$
By Schur's theorem mentioned in the beginning, $[F,N]$ is contained
in $S$ and it is straight forward to see that $[F,N]^n$ is contained
in $\myKernel(F,n)\cap S^n$ (both regarded as subgroups of $F^n$).
Hence this conjugate lies in $H$ and from this calculation it is easy to
deduce the normality of $H$ in $F^n$.

We already mentioned that $\myKernel(F,n)\cap S^n$ is a normal 
subgroup of $F^n$ contained in $\myKernel(F,n)$ and from the
presentation $G^\ast\iso F/S$ we easily get an isomorphism
\begin{equation}
\label{coverthmlabel}
\begin{array}{ccc}
   \myKernel(G^\ast,n)&\iso&\myKernel(F,n)\,/\,(\myKernel(F,n)\cap S^n)\,.
  \end{array}
\end{equation}
Since the image of $N$ in $F/S\iso G^\ast$ is $M$, we see that 
$\myKernel(N,n)$ maps to $\myKernel(M,n)$ inside
$\myKernel(G^\ast,n)$ under the isomorphism (\ref{coverthmlabel}).
Elements of the form $(1,...,m,1,...,m^{-1},1,...)$ with $m\in M$ 
generate $\myKernel(M,n)$, which is also true for $n=2$ since
$M$ is Abelian.
Such elements of $\myKernel(G^\ast,n)$ can be lifted to elements
of $\myKernel(F,n)$ lying inside $\myKernel(N,n)$.
Hence $\myKernel(N,n)$ maps surjectively onto $\myKernel(M,n)$.

Putting these results together we obtain an isomorphism
$$\begin{array}{ccc}
   \myKernel(G^\ast,n)/\myKernel(M,n) &\iso& \myKernel(F,n) / H \,.
  \end{array}
$$

By definition, the quotient of $\myKernel(F,n)$ by the normal closure
$\ll\myKernel(N,n)\gg$ of $\myKernel(N,n)$ inside $F^n$ is 
$\myBigKernel(G,n)$.
Clearly, $\myKernel(N,n)$ is contained in $H$ and since $H$ is normal
in $F^n$, also $\ll\myKernel(N,n)\gg$ is contained in $H$.
Hence the inclusion of these two normal subgroups induces a surjective
homomorphism
$$\begin{array}{ccccc}
   \psi & : & \myBigKernel(G,n) &\to& \myKernel(G^\ast,n)/\myKernel(M,n)\,.
  \end{array}
$$

On the other hand, both groups are extensions of
$H_2(G,\ZZ)$ by $\myKernel(G,n)$.
The surjective map $\psi$ induces an isomorphism between the 
kernel and the quotient of these extensions.
Hence $\psi$ is an isomorphism.
\qed

\section{Abelian groups }
\label{abeliansection}

We will now see that already
in the case of finite Abelian groups it is quite
difficult to determine the structure of $\myBigKernel(G,n)$.
Thanks to Proposition \ref{mykernelAbelian}, we know 
that $\myKernel(G,n)$ is isomorphic to $G^{n-1}$.

The following proposition implies that we may restrict ourselves
to $p$-groups: 
\begin{Proposition}
  \label{sylow}
  Let $G$ be a finite nilpotent group and $n\geq2$.
  Let $S_p$ be its unique Sylow $p$-subgroup.
  There exists an isomorphism
  $$\begin{array}{ccc}
    \myBigKernel(G,\,n)&\iso&\prod_{p}\,\myBigKernel(S_p,\,n),
  \end{array}$$
  where $p$ runs over all prime numbers.
  More precisely, the short exact sequence (\ref{mybigkerneleq})
  for $G$ is the product of the 
  short exact sequences (\ref{mybigkerneleq}) taken over all its
  Sylow $p$-subgroups $S_p$.
\end{Proposition}

\prf
By functoriality, there exists a commutative diagram with
exact rows
$$\begin{array}{cccccccccl}
  1&\to&\prod_p\,H_2(S_p)&\to&\prod_p \myBigKernel(S_p,\,n)&\to&
  \prod_p \myKernel(S_p,\,n)&\to&1\\
  &&\downarrow{\scriptstyle\varphi_1}
  &&\downarrow{\scriptstyle\varphi_2}
  &&\downarrow{\scriptstyle\varphi_3}\\
  1&\to&H_2(G)&\to&\myBigKernel(G,\,n)&\to&
  \myKernel(G,\,n)&\to&1&.\\
\end{array}$$
From \cite[Corollary 2.2.11]{ka} and Proposition \ref{mykernelproperties}
it follows that $\varphi_1$ and $\varphi_3$ are isomorphisms.
Hence $\varphi_2$ is an isomorphism.
\qed\medskip

First, we deal with $n=2$, which is the most interesting case from the
point of view of group theory.
Since $G$ is Abelian, $\myKernel(G,2)$ is isomorphic to $G$ and 
(\ref{mybigkerneleq}) becomes
\begin{equation}
\label{abelianeq}
\begin{array}{ccccccccc}
 0&\to& H_2(G,\,\ZZ) &\to& \myBigKernel(G,2) &\to& G &\to&1
\end{array}
\end{equation}

\begin{Proposition}
  \label{Abeliann=2}
  Let $G$ be a finite Abelian $p$-group.
  \begin{enumerate}
   \item If $p\neq2$ then $\myBigKernel(G,2)$ is a cover
  of $G$ via the short exact sequence (\ref{abelianeq}).
  \item If $G$ is an elementary Abelian $2$-group then
  $\myBigKernel(G,2)$ is an elementary Abelian $2$-group.
  More precisely, it is the product of $G$ and $H_2(G,\ZZ)$.
  \end{enumerate}
\end{Proposition}

\prf
Let $G$ be an Abelian $p$-group with $p\neq2$.
Every cover of $G$ is nilpotent of class at most $2$.
The same is true for $\myBigKernel(G,2)$ by (\ref{abelianeq}). 
Hence for arbitrary elements in these groups the commutator
relation $[x^i,y^j]=[x,y]^{ij}$ holds true.

To prove our statement we have to show that $H_2(G,\ZZ)$ in
(\ref{abelianeq}) lies in the commutator subgroup of 
$\myBigKernel(G,2)$.
We choose an arbitrary cover $G^\ast\,\to\,G$ with 
kernel $M$ and exhibit 
$\myBigKernel(G,2)$ as in Theorem \ref{coverthm}.
It is enough to prove that 
$\langle (m,1)\,|\,m\in M\rangle$ lies inside the commutator
subgroup of $\myKernel(G^\ast,2)$ modulo elements of
$\myKernel(M,2)$

Since $G$ is Abelian, $M$ coincides with the commutator subgroup
of $G^\ast$ and so we have to check that the commutator subgroup
of $G^\ast$ is a subgroup of the commutator subgroup of 
$\myKernel(G^\ast,2)$ modulo elements of $\myKernel(M,2)$.
Given $x,y\in G^\ast$, the elements $(x,x^{-1})$ and $(y,y^{-1})$
lie in $\myKernel(G^\ast,2)$ and hence
$([x,y],[x^{-1},y^{-1}])$ lies inside the commutator subgroup
of $\myKernel(G^\ast,2)$.
Modulo $\myKernel(M,2)$ this element is congruent to
$([x,y]^2,1)$.

Since $G$ has odd order, also $M$ has odd order by Schur's theorem,
cf. \cite[Theorem 2.1.5]{ka}.
Hence $([x,y],1)$ is a power of $([x,y]^2,1)$.
Thus, modulo elements of $\myKernel(M,2)$, the element $([x,y],1)$ 
lies in the commutator subgroup of $\myKernel(G^\ast,2)$.

Now, let $G$ be an elementary Abelian $2$-group.
To prove the remaining statement we can either proceed as above
or we copy the first part of the proof of 
Proposition \ref{centre} below.
%
\qed\medskip

\begin{Definition}
  Let $G$ be an Abelian group.
  We say that $G$ has a {\em natural cover} if
  $\myBigKernel(G,2)$ is a cover of $G$. 
  In this case, we will also refer to 
  $\myBigKernel(G,2)$ as {\em the natural cover of $G$}.
\end{Definition}

The following result is an immediate corollary of the previous proposition.

\begin{Theorem}
  A finite Abelian group of odd order possesses a natural cover.
  
  More precisely, let $G^\ast$ be an arbitrary cover of an Abelian
  group $G$ of odd order and let $M$ be the kernel of $G^\ast\,\to\,G$.
  Then we obtain the natural cover of $G$ as a subquotient of
  $(G^\ast)^2$ via
  $$
  \myBigKernel(G,2)\,\iso\,
  \langle (g,g^{-1}), (m,1) \,|\, g\in G^\ast, m\in M \rangle /
  \langle (m,m^{-1})\,|\,m\in M\rangle \,.
  $$
  Thus, the natural cover can be obtained from an arbitrary cover $G^\ast$
  by a kind of warped Baer sum 
  of $G^\ast$ with itself.
\end{Theorem}

\prf
By Proposition \ref{sylow} and Proposition \ref{Abeliann=2}
Abelian groups of odd order have natural covers.
The definition of $\myKernel(G^\ast,2)$ and Theorem \ref{coverthm}
give the explicit construction of the natural cover starting from
an arbitrary one.
\qed\medskip

We denote the centre of a group $G$ by $Z(G)$.
We denote its Frattini subgroup, i.e. the intersection of all maximal
subgroups of $G$, by $\Phi(G)$.

We recall that a $p$-group is called {\em special} if its centre is
equal to its commutator and its Frattini subgroup.
A special $p$-group is called {\em extra-special} if its centre
is cyclic.

\begin{Proposition}
Let $p$ be an odd prime number.
\begin{enumerate}
\item The natural cover of $\Cycl{p}$ is just $\Cycl{p}$ itself.
\item The natural cover of $(\Cycl{p})^2$ is the extra-special
      group of order $p^{1+2}$ and exponent $p$.
\end{enumerate}
\end{Proposition}

\prf
The first statement follows from Corollary \ref{mybigkernelAbelian}.

The Schur multiplier of $G=(\Cycl{p})^2$ is $\Cycl{p}$.
Hence $\myBigKernel(G,2)$ is a non-Abelian group of order $p^3$.
Such a group is necessarily extra-special.
The unique extra-special group $G^\ast$ of order $p^{1+2}$ and exponent $p$
is a cover of $G$.
Since $\myBigKernel(G,2)$ is a quotient of $\myKernel(G^\ast,2)$ by
Theorem \ref{coverthm}, the group $\myBigKernel(G,2)$ has exponent $p$.
This is enough to identify $\myBigKernel(G,2)$ as the 
unique extra-special group of order $p^{1+2}$ and exponent $p$.
\qed\medskip

For applications to algebraic geometry, especially the case $n\geq3$ is relevant.
Using Proposition \ref{mykernelAbelian}, the extension (\ref{mybigkerneleq}) 
becomes
\begin{equation}
\label{abelianbigeq}
\begin{array}{ccccccccc}
 0&\to& H_2(G,\,\ZZ) &\to& \myBigKernel(G,n) &\to& G^{n-1} &\to&1
\end{array}
\end{equation}
Since $n\geq3$, the group $H_2(G,\,\ZZ)$
lies inside the commutator subgroup of $\myBigKernel(G,n)$ by
Theorem \ref{mybigkernelthm}.

\begin{Proposition}
  Let $G$ be an Abelian $p$-group and $n\geq3$ a natural number.
  \begin{enumerate}
  \item Unless $G$ is cyclic, the group $\myBigKernel(G,n)$ is nilpotent
    of class $2$.
  \item The commutator subgroup of $\myBigKernel(G,n)$ is equal to 
    $H_2(G,\,\ZZ)$ embedded via (\ref{abelianbigeq}).
  \item The Frattini subgroup $\Phi(\myBigKernel(G,n))$ is an extension of
    $H_2(G,\,\ZZ)$ by $\Phi(G)^{n-1}$.
  \end{enumerate}
\end{Proposition}

\prf
If $G$ is not cyclic, then $H_2(G,\ZZ)$ does not vanish by Schur's
theorem, cf. \cite[Corollary 2.2.12]{ka}.
Since $\myBigKernel(G,n)$ is a central extension of two Abelian groups,
it is nilpotent of class at most $2$. 
Since $H_2(G,\ZZ)$ lies inside the commutator subgroup of $\myBigKernel(G,n)$
by Theorem \ref{mybigkernelthm},
the group $\myBigKernel(G,n)$ is not Abelian.

As $n\geq3$, the group $H_2(G,\ZZ)$ lies inside the commutator 
subgroup, which shows one inclusion.
On the other hand, the quotient of $\myBigKernel(G,n)$ by $H_2(G,\ZZ)$
is Abelian, showing the other inclusion.

Clearly, $\Phi(\myBigKernel(G,n))$ maps onto 
$\Phi(G^{n-1})$ via (\ref{abelianbigeq}).
It is a general fact that $I:=Z(G^\ast)\cap [G^\ast,G^\ast]$ is contained in
$\Phi(G^\ast)$.
Since $H_2(G,\ZZ)$ is contained in $I$ (in fact, they are equal in our
case), it follows
that $\Phi(\myBigKernel(G,n))$ is an extension of
$H_2(G,\ZZ)$ by $\Phi(G^{n-1})\iso \Phi(G)^{n-1}$.
\qed\medskip

The structure of the centre of $\myBigKernel(G,n)$ is much trickier.
In fact, it depends on $n$.

\begin{Proposition}
 \label{centre}
 Let $G$ be an Abelian $p$-group and $n\geq3$ a natural number.
 Let $G^\ast\to G$ be an arbitrary cover of $G$ and denote by $Z$ the
 image of the centre $Z(G^\ast)$ inside $G$.
 \begin{enumerate}
  \item If the exponent of $G$ divides $n$, then $\myBigKernel(G,n)$ is
    the direct product of $G$ and $\myBigKernel(G,n-1)$.
  \item If $p$ does not divide $n$, then the centre of $\myBigKernel(G,n)$
    is isomorphic to the product of $H_2(G,\ZZ)$ and $\myKernel(Z,n)$.
  \end{enumerate}
\end{Proposition}

\prf
We choose a cover $G^\ast\,\to\,G$ with kernel $M$ and define 
$Z$ as in the statement of the proposition.
By Theorem \ref{coverthm}, the quotient of
$\myKernel(G^\ast,n)$ by $\myKernel(M,n)$ is isomorphic to
$\myBigKernel(G,n)$.

Suppose that the exponent of $G$ divides $n$.
Then we obtain a well-defined injective homomorphism $\Delta$ 
from $G^\ast$ to $\myKernel(G^\ast,n)$ that sends
$g$ to $(g,....,g)$.
Since also $M$ has an exponent which divides $n$,
we conclude that the intersection 
$\Delta(M)\cap\myKernel(M,n)$ is equal to $\Delta(M)$.
Hence we obtain $G^\ast/M= G$ as a central subgroup of
$\myBigKernel(G,n)$.
This subgroup maps to a diagonally embedded $G$ inside
$\myKernel(G,n)$ (under the map (\ref{mybigkerneleq}))
and hence we can split the induced
injective map from $G$ to $\myBigKernel(G,n)$.
Thus, $G$ is a direct factor of $\myBigKernel(G,n)$ and
it is easy to see that the quotient is in fact
isomorphic to $\myBigKernel(G,n-1)$.
This also works for $n=2$, but then the quotient of
$\myBigKernel(G,2)$ by $G$ is equal to $M$.

To prove the second assertion we now assume that $p$ does not
divide $n$.
The preimage of the centre of $\myBigKernel(G,n)$ in $\myKernel(G^\ast,n)$
consists of those elements of $\myKernel(G^\ast,n)$ for which every
commutator lies in $\myKernel(M,n)$.
Hence this preimage is equal to 
$$
\tilde{Z} \,:=\, \left\{ 
(h_1,...,h_n)\in\myKernel(G^\ast,n)\,|\, \sum_{i=1}^n [g_i,h_i]=0 \,\,
\forall (g_1,...,g_n)\in\myKernel(G^\ast,n)
\right\}
$$
Let $(h_1,...,h_n)$ be an element of $\tilde{Z}$.
For $g\in G^\ast$, the element $(1,...,g,1,...,g^{-1},1,...)$ lies in 
$\myKernel(G^\ast,n)$ and
we obtain $[h_i,g]=[h_j,g]$ for all $i,j$.
In particular, if $h_i\in Z(G^\ast)$ for some $i$, then 
$h_i\in Z(G^\ast)$ for all $i$.

Assume there exist an element 
$\vec{h}:=(h_1,...,h_n)$ of $\tilde{Z}$ 
with $h_i\not\in Z(G^\ast)$ for some $i$, say $i=1$.
Inside $(G^\ast)^n$ we can write this element as product of
$(h_1,...,h_1)$ by $\vec{z}:=(1,h_1^{-1}h_2,...,h_1^{-1}h_n)$,
where all entries of $\vec{z}$ lie in $Z(G^\ast)$.
The sum over all components of $\vec{z}$ is an element 
$z'$ of $Z(G^\ast)$.
Since $p$ does not divide $n$, there exists a power $z''$ of
$z'$ that is an $n$.th root of $z'$.
We define $\vec{h}':=(h_1\cdot z'',...,h_1\cdot z'')$
and $\vec{z}':=(z'',...,z'')^{-1}\vec{z}$.
Then $\vec{h}=\vec{h}'\cdot\vec{z}'$.
We arranged $\vec{z}'$ in such a way that it lies
in $\myKernel(G^\ast,n)$ and so also $\vec{h}'$ lies in
$\myKernel(G^\ast,n)$.
Every component of $\vec{h}'$ is equal to $h_1 z''$.
Since $h_1 z''$ does not lie $[G^\ast,G^\ast]$ 
(this group is contained in the centre of $G^\ast$ and we
assumed $h_1\not\in Z(G^\ast)$),
also $(h_1 z'')^n$ does not lie in $[G^\ast,G^\ast]$ (using
again the fact that $n$ is coprime to $p$).
Hence $\vec{h}$ does not lie in $\myKernel(G^\ast,n)$,
a contradiction.

We conclude that the centre of $\myBigKernel(G,n)$
is the image of $Z(G^\ast)^n\cap\myKernel(G^\ast,n)$, i.e.
we have to compute its quotient by $\myKernel(M,n)$. 
This group, however, is easily seen to be an extension
of $M$ by $\myKernel(Z,n)$.
\qed

\begin{Corollary}[Read's theorem for Abelian groups]
  Let $G^\ast\,\to\,G$ be a cover of an Abelian group $G$.
  Then the image $Z$ of the centre $Z(G^\ast)$ inside $G$ does not depend on
  the choice of the cover $G^\ast$.
\end{Corollary}

\prf
For a natural number $n\geq3$ that is coprime to the order of $G$,
the centre of $\myBigKernel(G,n)$ is a product of $H_2(G,\ZZ)$ and
$\myKernel(Z,n)$.
Since this centre does not depend on the choice of a cover, also
$Z$ is independent of it.
\qed

\begin{Remark}
  In view of the last corollary it seems natural to ask whether
  $\myBigKernel(G,n)$ captures interesting data about all cover
  groups of a given (not necessarily Abelian) group $G$.
  For example, Schur's theorem that $[G^\ast, G^\ast]$ is an invariant
  of $G$ and does not depend on the choice of the cover $G^\ast$ 
  also follows quite easily from Theorem \ref{coverthm}.
\end{Remark}
  
\section{Fundamental groups}
\label{surfacesection}

We now sketch how $\myKernel(G,n)$ and $\myBigKernel(G,n)$ are connected
to fundamental groups of algebraic surfaces and Moishezon's
programme to find new invariants for algebraic surfaces.
For details and references we refer to \cite{li}.

For complex curves, it is already known since the 19th century, that
their fundamental groups classify them up to diffeomorphism.
However, although the Italian school classified complex algebraic 
surfaces of special type in the early 20th century, 
not much is known about surfaces of general type.

There is an approach towards a finer classification 
that uses embeddings of surfaces into large projective spaces:
Let $X$ be a smooth projective surface and $\mathcal L$
a sufficiently ample line bundle on $X$.
Then we embed $X$ via $\mathcal L$ into some projective space $\PP^N$.
After that we choose a generic codimension three linear subspace in
this $\PP^N$ and consider the projection $\pi$ away from this space.
This is a rational map from $\PP^N$ onto $\PP^2$.
The restriction 
$$\begin{array}{ccccc}
f\,:=\,\pi|_X&:&X&\to&\PP^2
\end{array}$$
is a finite map, called a {\em generic projection}.

We denote by $n$ the degree of $f$ and by $D$ its branch locus.
If we know the fundamental group $\pi_1(\PP^2-D)$ and the monodromy 
morphism $\psi:\pi_1(\PP^2-D)\,\to\,\Sym_n$, where $\Sym_n$ 
denotes the symmetric group on $n$ letters, we can reconstruct $X$.
Thus, if we could extract invariants from these fundamental groups
we would get a much finer classification of algebraic surfaces.
However, these groups $\pi_1(\PP^2-D)$ are huge and may have a
rather complicated structure, although a conjecture of Teicher
states that they are almost-solvable.

Also this was known for some time, but it could not be
used effectively since it was too difficult to compute these 
fundamental groups.
However, the braid group techniques introduced by Moishezon in \cite{mo}, 
and refined later on by Teicher and others, made it possible to 
compute these groups $\pi_1(\PP^2-D)$ in many cases.
\medskip

One such invariant for $X$ depending on the choice of
an embedding of $X$ into projective space is the fundamental group
of the Galois closure of this generic projection:
To a generic projection $f$ of degree $n$ we associate 
its so-called {\em Galois closure}
$$\begin{array}{ccccc}
\Xgal &:=& \overline{\{ (x_1,...,x_n)\,\,|\,\,x_i\neq x_j,\,f(x_i)=f(x_j) \}}&\subseteq&X^n .
\end{array}$$
This turns out to be a smooth projective surface. 
In most cases it is of general type.

Apart from their connection with generic projections, there
is another reason why Galois closures are interesting:
Using Galois closures of generic projections gives one of the few known ways to 
construct series of surfaces of general type with positive index, 
i.e. the Chern numbers of $\Xgal$ fulfil $c_1{}^2>2c_2$.
For some time it was believed that surfaces of general type with positive 
index should have infinite fundamental groups.
The first counter-examples to this conjecture were given by Moishezon and Teicher via
computing fundamental groups of Galois closures of generic projections from
$X=\PP^1\times\PP^1$, \cite{mote}.
\medskip

Hence, determining $\pi_1(\Xgal)$ is interesting from the point of view
of fundamental groups of surfaces of general type.
Also, these groups should give new invariants of $X$ as they occur
as certain naturally defined subquotients of $\pi_1(\PP^2-D)$.
In \cite{li} we partly simplified the calculations of \cite{mote} 
which led to $\myBigKernel(G,n)$:

We fix a generic projection $f:X\to\PP^2$ of degree $n$ and denote by $\Xgal$ its
associated Galois closure $\fgal:\Xgal\to\PP^2$.
We fix a generic line in $\PP^2$ and denote its complement by $\Aff^2$.
Then we denote the inverse images of $f^{-1}(\Aff^2)$ and $\fgal^{-1}(\Aff^2)$ by
$\Xaff$ and $\Xaffgal$, respectively.
Since $f$ is generic the Galois group of $\Xgal$ over $\PP^2$ is the whole 
symmetric group $\Sym_n$.
This group acts on $\Xgal$ and $\Xaffgal$ and we can form the quotient
$$\begin{array}{lcl}
  \blowup\Xgal\,/\,\Sym_{n-1} &\iso& X\\
  \Xaffgal\,/\,\Sym_{n-1} &\iso& \Xaff\,.
\end{array}$$
There are $n$ distinct embeddings of $\Sym_{n-1}$ into $\Sym_n$ yielding $n$ 
distinct isomorphisms and $n$ distinct induced maps on fundamental groups, all of
which are surjective:
$$\begin{array}{lcl}
  \blowup\pi_1(\Xgal) &\to& \pi_1(X)\\
  \pi_1(\Xaffgal) &\to& \pi_1(\Xaff)\,.
\end{array}$$
Combining these $n$ homomorphisms, we obtain a map from $\pi_1(\Xgal)$ to
$\pi_1(X)^n$, and similarly for $\pi_1(\Xaffgal)$.
The following result determines the images of these maps.
\begin{Theorem}
\label{thm1}
There exist surjective homomorphisms
$$\begin{array}{lcl}
  \blowup\pi_1(\Xgal) &\to& \myKernel(\pi_1(X),\,n)\\
  \pi_1(\Xaffgal) &\to& \myKernel(\pi_1(\Xaff),\,n)\,.
\end{array}$$
\end{Theorem}
The arguments in the proof of Theorem \ref{thm1} can be formalised in such 
a way that the result
remains true for \'etale fundamental groups and generic projections 
defined over arbitrary algebraically closed fields of 
characteristic $\neq2,3$.
Of course, one has to modify the statement for
$\pi_1(\Xaffgal)$ over fields of positive characteristic a little bit
since the affine plane is then no longer simply-connected.

Over the complex numbers the algorithm of Zariski and van~Kampen provides us
with a presentation of fundamental groups of complements of curves in $\Aff^2$
or $\PP^2$.
Applying it to the branch curve $D$  of $f$ we find the fundamental group
$\pi_1(\Xaffgal)$ as a subquotient of $\pi_1(\Aff^2-D)$.
Combining this presentation with Theorem \ref{thm1} we obtain
\begin{Theorem}
\label{phdthm}
There exists a surjective homomorphism
$$\begin{array}{lcl}
\pi_1(\Xaffgal) &\to& \myBigKernel(\pi_1(\Xaff),\,n)\,.
\end{array}$$
The group $\pi_1(\Xgal)$ is a quotient of $\pi_1(\Xaffgal)$ by a cyclic central subgroup.
\end{Theorem}

In all known examples, where the generic projection was defined via a sufficiently ample
line bundle, the map of Theorem \ref{phdthm} is in fact an isomorphism.

This suggests to use Galois closures of generic projections to construct algebraic 
surfaces with interesting fundamental groups.
For example, starting from a surface with Abelian fundamental group,
iterated Galois closures should produce surfaces with nilpotent fundamental groups
of large class. 
Another project would be to obtain new surfaces with fundamental groups that
are not residually finite.

Whether the map of Theorem \ref{phdthm} is an isomorphism in all cases or 
at least for large class of surfaces or generic projections is not clear at the moment,
although this is true in all known examples.
In any case, we have a quotient of the group we are interested in
and the appearance of covering groups in connection with these fundamental groups
is quite surprising.
Therefore, it is indispensable to have a better understanding
of $\myBigKernel(-,n)$ to comprehend these fundamental groups.


\begin{thebibliography}{XXXX}
  \bibitem[Br]{br} K.~S.~Brown, {\em Cohomology of Groups},
    GTM 87, Springer (1982).
  \bibitem[Ka]{ka} G.~Karpilovsky, {\em The Schur Multiplier},
    LMS Monographs 2, Oxford University Press (1987).
  \bibitem[Li]{li} C.~Liedtke, {\it On Fundamental Groups of
    Galois Closures of Generic Projections}, 
    Bonner Mathematische Schriften  367 (2004).
  \bibitem[Mo]{mo} B.~Moishezon, {\it Stable branch curves and braid monodromies},
    Springer LNM 862 (1981), 107-192.
  \bibitem[MoTe]{mote} B.~Moishezon, M.~Teicher, {\it Simply-connected algebraic
    surfaces of positive index},
    Invent. Math. 89, 601-643 (1987).
\end{thebibliography}
\end{document}